%% file: Control_Class.tex
\documentclass[final]{siamltex}

\usepackage{amsmath} 
\usepackage{amsfonts}

\usepackage{float}
\usepackage{graphicx}
\usepackage{subcaption}
\usepackage{url}
\usepackage{color}
\usepackage{bm}
\usepackage{multirow}

\usepackage{listings}
\usepackage{color}

\definecolor{mygreen}{rgb}{0,0.6,0}
\definecolor{mygray}{rgb}{0.5,0.5,0.5}
\definecolor{mymauve}{rgb}{0.58,0,0.82}

\usepackage{hyperref}
\hypersetup{
    colorlinks,
    linkcolor={red},
    citecolor={blue},
    urlcolor={blue}
}

\title{Automatic Differentiation for All-at-once Systems Arising in Certain PDE-Constrained Optimization Problems}

\author{Santolo Leveque\thanks{CRM Ennio De Giorgi, Scuola Normale Superiore, Piazza dei Cavalieri 3, 56126, Pisa, Italy ({\tt santolo.leveque@sns.it})}
\and
James R. Maddison\thanks{School of Mathematics and Maxwell Institute for Mathematical Sciences, The University of Edinburgh, James Clerk Maxwell Building, The King's Buildings, Peter Guthrie Tait Road, Edinburgh, EH9 3FD, United Kingdom (\tt j.r.maddison@ed.ac.uk)}
\and
John W. Pearson\thanks{School of Mathematics and Maxwell Institute for Mathematical Sciences, The University of Edinburgh, James Clerk Maxwell Building, The King's Buildings, Peter Guthrie Tait Road, Edinburgh, EH9 3FD, United Kingdom ({\tt j.pearson@ed.ac.uk})}}

\begin{document}
\maketitle

\begin{abstract}
An automated framework is presented for the numerical solution of optimal
	control problems with PDEs as constraints, in both the stationary and
	instationary settings. The associated code can solve both linear and
	non-linear problems, and examples for incompressible flow equations
	are considered. The software, which is based on a Python interface to
	the Firedrake system, allows for a compact definition of the problem
	considered by providing a few lines of code in a high-level language.
	The software is provided with efficient iterative linear solvers for
	optimal control problems with PDEs as constraints. The use of advanced
	preconditioning techniques results in a significant speed-up of the solution
	process for large-scale problems. We present numerical examples of the
	applicability of the software on classical control problems with PDEs as
	constraints.
	
\end{abstract}



\pagestyle{myheadings}
\thispagestyle{plain}
\markboth{S. LEVEQUE, J. R. MADDISON, AND J. W. PEARSON}{AD FOR PDE-CONSTRAINED OPTIMIZATION}

\section{Introduction}\label{sec_1}
Optimal control problems with partial differential equations (PDEs) as
	constraints arise very naturally in industrial and real-life applications.
%	A classical example is given by the heating of a tissue during a medical
%	treatment: one wishes to drive the temperature of the body towards a
%	desired one, with a minimal cost in terms of energy.
	Such problems are formulated as the minimization of
	a cost functional subject to the system of PDEs that describe the
	physics considered. Problems belonging to this class have been applied
	to shape optimization \cite{Abraham_Behr_Heinkenschloss, Orozco_Ghattas,
	Quarteroni_Rozza}, medical imaging and tomography \cite{Arridge,
	Cheney_Isaacson_Newell, Klibanov_Lucas}, mathematical finance
	\cite{Bouchouev_Isakov, Egger_Engl}, reaction--diffusion control
	for chemical processes \cite{Barthel_John_Troltzsch, Griesse_Volkwein},
	optimal design of semiconductor devices \cite{Hinze_Pinnau},
	and flow control in porous media \cite{Ewing_Lin, Hasan_Foss_Sagatun,
	Sudaryanto_Yortsos}.

Since the early work by Lions \cite{Lions}, many researchers have devoted
	effort to devising methods for the numerical solution of the
	optimal control of PDEs. Frequently, the most complex task when solving
	this type of problem is to devise efficient algorithms for the solution of
	the first-order optimality conditions (also known as Karush--Kuhn--Tucker,
	or simply KKT, conditions), which consist of a system of discretized PDEs
	that the solution has to satisfy. The KKT conditions can be derived by
	employing two approaches. On one hand, one can consider the problem as 
	an optimization problem in an infinite-dimensional Hilbert space, and
	derive KKT conditions by employing classical functional analysis; then,
	the KKT conditions are discretized, obtaining a sequence of linear systems
	to be solved at each non-linear iteration. This approach is called
	\emph{optimize-then-discretize}. On the other hand, one can first discretize
	the cost functional and the system of PDEs considered, and then derive
	the corresponding KKT conditions of the resulting finite-dimensional
	optimization problem. This approach is called
	\emph{discretize-then-optimize}. Either way, one has to solve at each
	non-linear iteration a very large linear system. In fact, the dimension
	of the system to be solved could become prohibitive when one adopts
	the so called \emph{all-at-once} approach. In this case, one solves
	the KKT conditions for all the solutions concurrently. In particular, when
	considering a time-dependent PDE as a constraint, the discretization
	of the differential operator has to take into account the time
	derivative, thus even storing the right-hand side of the linear system
	to be solved may become unfeasible.

Due to the high dimensionality of the KKT conditions, a common approach
	is employing a gradient-type method, in which simpler (uncoupled forward
	and adjoint) problems are solved until an accurate enough solution is
	found; see, among many examples, \cite{Geiersbach_Wollner, Gotschel_Minion}.
	For time-dependent initial-value problems, gradient-type methods solve
	two instationary PDEs, one marching forward in time, the other backward,
	see, for instance, \cite{Gotschel_Minion}. For this class of problems,
	we also point out multiple shooting methods \cite{Hein05}, parareal schemes
	\cite{MaTu02, MSS10}, and ideas from instantaneous control
	\cite{CHK99, Hinz00}, among many other effective approaches.

Although one has to solve smaller problems, gradient descent may be expected
	to converge at first order for sufficiently smooth functions, meaning
	that one requires many iterations in order to obtain a very accurate
	solution. In this sense, one would favour adopting Newton-type methods,
	that are able to attain at least superlinear convergence, see, for instance,
	\cite{Bergounioux_Haddou_Hintermuller_Kunisch, Hintermuller_Ito_Kunisch,
	Hintermuller_Hinze, Hinze_Koster_Turek, Ito_Kunisch_2003, Schiela_Weiser,
	Ulbrich_Ulbrich, Weiser, Weiser_Schiela}. Further, some Newton-type methods
	are also able to deal with issues of the regularity of the solution of an
	optimal control problem with PDEs as constraints and additional algebraic
	constraints on the state variable, see, for instance,
	\cite{Hintermuller_Ito_Kunisch, Ito_Kunisch_2003}. In fact, for this
	class of problems one can prove that the Lagrange multiplier corresponding
	to the bounds on the state variable is only a Borel measure,
	see \cite{Casas}.
%In order to deal with this issue, a regularization technique is employed within the optimization framework, obtaining in this way a semismooth Newton method, see, for instance, \cite{Hintermuller_Ito_Kunisch, Ito_Kunisch_2003}.

Despite the great effort being devoted to the construction of methods for
	the solution of optimal control with PDEs as constraints, the implementation
	of flexible open-source codes for the solution of this class of problems
	forms a pressing challenge within the community. Regarding this, we recall
	here the work \cite{Funke_Farrell}, in which the authors present a framework
	based on automatic differentiation (AD) for the solution of this
	class of problems. AD represents a powerful tool for the solution
	of problems involving PDEs by finite elements: by employing a
	domain-specific language, high-level AD exploits the variational form of the
	problem in order to simplify the discretization of the problem
	\cite{Farrell_Ham_Funke_Rognes}. In this way, the user is required to
	provide only a few lines of code of a high-level language. The work
	we present here points toward this direction: by employing AD tools,
	we develop software that allows the user to define the optimal
	control problem in few lines of code. We would like to note that our work
	differs from that of \cite{Funke_Farrell}. In fact, the software
	implemented in \cite{Funke_Farrell} allows the user to impose additional
	algebraic constraints on the variables, while our software is currently
	able to solve only problems without additional algebraic constraints.
	However, as opposed to the linear solver employed in \cite{Funke_Farrell},
	which is based on a block-LU decomposition, our software employs
	advanced preconditioning techniques in order to solve for the discretized
	KKT conditions. Recently, advanced preconditioners for time-dependent
	all-at-once optimization problems have been developed in
	\cite{Leveque_Pearson_2021, Leveque_Pearson_2022, Pearson_Stoll_Wathen}.
	These preconditioners rely on the ability to apply a block bi-diagonal or
	a block lower-triangular preconditioner to the systems arising from
	backward Euler and trapezoidal rule (Crank--Nicolson type) discretizations.
	A key development described in this article is the ability to apply
	these preconditioners to more general time-dependent problems. In
	particular, in this work the key details of the time-dependent
	PDE are described at a high-level, and automatic differentiation is
	applied to form the KKT system automatically. The preconditioners can then
	be applied automatically, in a similar way. This opens up these advanced
	preconditioners for application to a wider range of problems.

%We would like to note that the complex structure of the all-at-once linear
%	system arising when employing, for example, a trapezoidal rule in
%	time together with the corresponding optimal preconditioner results {\color{blue}in
%	a complicated classical AD code--the bi-diagonal transformations employed in
%	the preconditioner derived in \cite{Leveque_Pearson_2021} cannot be set implicitly in the
%	bilinear form passed to Firedrake}.
We would like to note that it is difficult to combine an optimal
	preconditioner for the KKT system of the problems considered here
	with a classical AD approach. For example, the optimal preconditioner
	derived in \cite{Leveque_Pearson_2021} employs a linear operator that can be
	factorized as the product of block bi-diagonal matrices as an approximation
	of the (negative) Schur complement. The implementation of such an
	approximation in a classical AD framework would result in a long
	code with a complicated structure. By contrast, the use of a Python
	interface allows our code to define the problem, derive the KKT conditions
	and the discretized system, and finally apply the preconditioner in a much
	simpler way. By keeping as much as possible in the Python layer, it is more
	straightforward to define new and more complicated preconditioners which
	make more detailed use of the structure of the block system arising in
	general all-at-once optimization problems.
	
The preconditioner developed in our software is based on the ``matching
	strategy''. Firstly employed for solving Poisson control problems in
	\cite{Pearson_Wathen_2012}, the matching strategy has been successfully
	applied to the solution of stationary convection-diffusion control problems
	in \cite{Pearson_Wathen_2013}, heat control problems in
	\cite{Leveque_Pearson_2021, Pearson_Stoll_Wathen}, and incompressible fluid
	flow control problems in \cite{Leveque_Pearson_2022}.

The software is available at \url{https://github.com/sleveque/control}.

This paper is structured as follows. In Section \ref{sec_2}, we introduce the
	problem we consider, that is, the optimal control of (stationary or
	instationary) PDEs. In Section \ref{sec_3}, we give details of our software,
	describing the preconditioner employed in the linear solver. In Section
	\ref{sec_4}, we present numerical examples of the type of problems that
	our software can handle, showing the lines of code that the user is
	required to provide. Finally, we give conclusions and avenues for future
	work in Section \ref{sec_5}.

\section{Problem Formulation}\label{sec_2}
In this section, we introduce the class of problems that our software can
	handle. For simplicity, we consider here only linear problems,
	noting that the reasoning may be applied to some linearizations of non-linear
	problems as well. The non-linear iterations that can be employed within our
	software are specified in Section \ref{sec_4}.

Given a suitable domain $\Omega \subset \mathbb{R}^d$, with $d\in \mathbb{N}$, and
	$\beta>0$, we consider the distributed optimal control problem of the
	following stationary PDE:
	\begin{equation}\label{cost_functional_stationary}
		\min_{v \in V, u \in U} J(v,u) := 
			\frac{1}{2}\| v - v_d \|^2_{L^2(\Omega)} +
			\frac{\beta}{2} \| u \|^2_{L^2(\Omega)}
	\end{equation}
	subject to
	\begin{equation}\label{stationary_PDE_CO}
		\mathcal{D}v = u + f \qquad \mathrm{in} \; \Omega,
	\end{equation}
	where $\mathcal{D}$ is a differential operator in a certain (given) space,
	equipped with suitable boundary conditions, and $V, U\subseteq L^2(\Omega)$
	are suitable Hilbert spaces. The function $v$ is called
	the state, and the function $u$ is called the control. Further, the function
	$v_d$ is called the desired state, while the term $f$ is a forcing
	function acting on the system.

Similarly, given also a final time $t_f>0$, we consider the distributed
	optimal control problem of instationary PDEs, defined as the minimization
	of the cost functional
	\begin{equation}\label{cost_functional_instationary}
		\min_{v \in V, u \in U} J_t(v,u) = \int_0^{t_f} J(v,u) \mathrm{d} t
	\end{equation}
	with $J(v,u)$ defined in \eqref{cost_functional_stationary}, subject to
	\begin{equation}\label{instationary_PDE_CO}
		\frac{\partial v}{\partial t}+\mathcal{D}v = u + f
			\qquad \mathrm{in} \; \Omega \times (0, t_f),
	\end{equation}
	equipped with some suitable initial and boundary conditions. We refer to
	\cite{Hinze_Pinnau_Ulbrich_Ulbrich, Lions, Troltzsch} for a detailed
	discussion of optimal control problems with PDEs as constraints.

We employ the discretize-then-optimize strategy when solving stationary
	control problems, while adopting a hybrid approach for the control of
	instationary PDEs. Specifically, if solving a stationary control problem
	we discretize the cost functional \eqref{cost_functional_stationary} and the
	constraints \eqref{stationary_PDE_CO}, and then derive the KKT conditions
	by employing classical constrained optimization theory. On the other hand,
	when solving an instationary control problem we employ the
	optimize-then-discretize strategy in the time variable, alongside a
	method of lines approach that utilizes suitable spatial discretizations
	of the terms in
	\eqref{cost_functional_instationary}--\eqref{instationary_PDE_CO}, thus
	obtaining a system of ordinary differential equations (ODEs) as constraints.
	We employ the finite element method for the discretization in space,
	adopting an all-at-once approach when deriving the discretized system.

We start by introducing the KKT conditions for the stationary case. Letting
	$\mathbf{M}_v^{-1}\mathbf{D}$ (for an appropriate scaling matrix
	$\mathbf{M}_v$) be a suitable discretization of $\mathcal{D}$,
	the discretization of
	\eqref{cost_functional_stationary}--\eqref{stationary_PDE_CO} reads as
	\begin{displaymath}
		\min_{\bm{v},\bm{u}} \bm{J}(\bm{v},\bm{u}) =
			\frac{1}{2}(\bm{v} - \bm{v}_d)^\top \mathbf{M}_1 (\bm{v} - \bm{v}_d)
			+ \frac{\beta}{2} \bm{u}^\top \mathbf{M}_2 \bm{u},
	\end{displaymath}
	subject to
	\begin{displaymath}
		\mathbf{D} \bm{v} = \mathbf{M}_v \bm{u} + \bm{f},
	\end{displaymath}
	for appropriate $\mathbf{M}_1$ and $\mathbf{M}_2$. Note that both
	$\mathbf{M}_1$ and $\mathbf{M}_2$ are symmetric positive definite.
	Note also that $\mathbf{M}_1=\mathbf{M}_2=\mathbf{M}_v$ if
	$v$ and $u$ are sought in the same finite element space,
	an assumption that we will make from now onwards\footnote{Although
	this assumption may seem limiting, it is essential from the point of view
	of the numerical linear algebra; in fact, to the best of the authors'
	knowledge, no optimal preconditioners are currently available for a
	more general setting.}. Here, $\bm{v}$ and
	$\bm{u}$ contain the numerical approximations of $v$ and
	$u$, respectively, while $\bm{f}$ is the discretization of the function
	$f$. Then, by introducing the adjoint variable
	$\bm{\zeta}$ and employing classical theory of constrained
	optimization, one can derive the following KKT conditions:
	\begin{equation}\label{discretized_KKT_stationary}
		\left\{
			\begin{array}{l}
				\mathbf{M}_v (\boldsymbol{v} - \boldsymbol{v}_d) +
					\mathbf{D}^* \boldsymbol{\zeta} = \boldsymbol{0},\\
				\mathbf{D} \boldsymbol{v} -
					\frac{1}{\beta} \mathbf{M}_v \boldsymbol{\zeta} -
					\boldsymbol{f} = \boldsymbol{0},
			\end{array}
		\right.
	\end{equation}
	with $\mathbf{D}^*$ arising from the adjoint term of $\mathcal{D}$.

We now introduce the KKT conditions for the instationary control problem. As
	mentioned, we first discretize the cost functional
	\eqref{cost_functional_instationary} and PDE
	\eqref{instationary_PDE_CO} in space. This results in
	\begin{displaymath}
		\min_{\bm{v},\bm{u}} \bm{J}_t(\bm{v},\bm{u}) = \int_0^{t_f}
			\left[ \frac{1}{2}(\bm{v} - \bm{v}_d)^\top \mathbf{M}_v (\bm{v}
			- \bm{v}_d) + \frac{\beta}{2} \bm{u}^\top \mathbf{M}_v
			\bm{u} \right] \mathrm{d}t,
	\end{displaymath}
	subject to
	\begin{equation}\label{discretized_constraints_instationary_control}
		\begin{array}{ll}
			\mathbf{M}_v \frac{\mathrm{d} \bm{v}}{\mathrm{d} t} +
				\mathbf{D} \bm{v} = \mathbf{M}_v \bm{u} + \bm{f}
				& \quad \mathrm{in} \; (0, t_f).
		\end{array}
	\end{equation}
	The constraints in \eqref{discretized_constraints_instationary_control} are
	a system of ODEs. In order to derive the first-order optimality conditions,
	we adopt an optimize-then-discretize strategy in the time variable.
	Introducing the adjoint variable $\bm{\zeta}$, by employing the Lagrangian
	method one can obtain the following KKT conditions:
	\begin{equation}\label{discretized_KKT_instationary}
		\left\{
			\begin{array}{ll}
				\mathbf{M}_v (\bm{v} - \bm{v}_d) - \mathbf{M}_v
					\frac{\mathrm{d} \bm{\zeta}}{\mathrm{d} t} +
					\mathbf{D}^* \bm{\zeta} = \bm{0}, &
					\quad \mathrm{in} \; (0, t_f),\\
				\mathbf{M}_v \frac{\mathrm{d} \bm{v}}{\mathrm{d} t} +
					\mathbf{D} \bm{v} - \frac{1}{\beta} \mathbf{M}_v
					\bm{\zeta} - \bm{f}
					= \bm{0}, & \quad \mathrm{in} \; (0, t_f).
			\end{array}
		\right.
	\end{equation}
	The last step consists of discretizing the time derivatives of $\bm{v}$
	and $\bm{\zeta}$. For time-dependent problems, the software allows the
	user to employ either a backward Euler or a trapezoidal rule discretization
	for the time derivatives.
%By doing so, we obtain the following discretized KKT
%	conditions:\color{blue}
%	\begin{equation}\label{discretized_KKT_instationary}
%		\left\{
%			\begin{array}{l}
%				\mathbf{M}_v (\bm{v} - \bm{v}_d) - \mathbf{M}_v
%					\mathcal{T}(\bm{\zeta})+
%					\mathbf{D}^* \bm{\zeta} = \bm{0},\\
%				\mathbf{M}_v \mathcal{T}(\bm{v}) + \mathbf{D} \bm{v}
%					- \frac{1}{\beta} \mathbf{M}_v \bm{\zeta} - \bm{f}
%					= \bm{0}.
%			\end{array}
%		\right.
%	\end{equation}\color{black}
%	Here, $\mathcal{T}(\cdot)$ is a suitable discretization of the time
%	derivative.

Finally, we write the all-at-once systems of the discretized KKT conditions
	given in \eqref{discretized_KKT_stationary} and
	\eqref{discretized_KKT_instationary}. The all-at-once system can be
	written as follows:
	\begin{equation}\label{general_block_system}
		\underbrace{\left[
			\begin{array}{cc}
				A & B_1^\top\\
				B_2 & - C
			\end{array}
		\right]}_{\mathcal{A}}
		\left[
			\begin{array}{c}
				\boldsymbol{v}\\
				\boldsymbol{\zeta}
			\end{array}
		\right]
		= \boldsymbol{b},
	\end{equation}
	for a suitable right-hand side $\boldsymbol{b}$. The $(1,1)$-block $A$
	is symmetric positive definite, while the $(2,2)$-block $-C$ is symmetric
	negative definite. Further, the blocks $B_1^\top$ and $B_2$ contain the
	discretization of the forward and the adjoint differential operators
	$\mathcal{D}$ and $\mathcal{D}^*$, respectively.
	Note that system \eqref{general_block_system} may form part of a larger
	system to be solved, as is the case when adding additional (differential--)
	algebraic constraints on the state variable, for instance. In fact, we also
	consider fluid flow control problems, such as Stokes or Navier--Stokes
	control problems, where additional incompressibility constraints are applied
	to the state variable, in addition to terms being applied to an additional
	pressure variable.

The automated differentiation tools of Firedrake allows the software to derive
	the adjoint discretized operator $\mathbf{D}^*$, asking the user to
	provide only the forward differential operator $\mathcal{D}$.

\section{The Software}\label{sec_3}
In this section, we give details of our software. We employ the Firedrake
	system
	\cite{Rathgeber_Ham_Mitchell_Lange_Luporini_Mcrae_Bercea_Markall_Kelley}
	to derive the finite element models of the problems considered,
	using the Python interface to PETSc \cite{Dalcin_Paz_Kler_Cosimo}
	for the derivation of the KKT conditions and the definition of the
	linear solvers.

\subsection{Firedrake}\label{sec_3_1}
Firedrake is an automated system employed for the numerical solution of PDEs.
	By employing the Unified Form Language (UFL)
	\cite{Alnaes_Logg_Olgaard_Rognes_Wells}, the user is able to express PDEs in
	weak form, which is then passed to a code generator which generates code
	which can, for example, assemble discretized operators. The symbolic
	representation of the differential operator through UFL allows the user
	to derive derivative information, including associated adjoint operators,
	in an automated way. From here, the user is required only to write few
	lines of code for the definition of the problem considered, by providing,
	aside from the weak formulation of the PDEs, the finite element space to
	which the solution belongs and the boundary conditions that the solution
	has to satisfy.

By making use of the Python interface, Firedrake is able to interact with
	PETSc, allowing the user to specify the linear and the non-linear iteration
	to employ for the solution of the discretized PDEs. The wide range of linear
	as well as non-linear solvers implemented in PETSc allows the user to build
	a problem-specific solver in few lines of code.

\subsection{The Control Class}\label{sec_3_2}
In this section, we discuss the extensions to the Firedrake system added in this
	work, showing how easily one can set up a control problem.

For simplicity, we first consider the following heat control problem:
	\begin{displaymath}
		\min_{v, u} \frac{1}{2} \int_0^{t_f} \| v - v_d \|^2_{L^2(\Omega)}
		\mathrm{d} t + \frac{\beta}{2} \int_0^{t_f} \| u \|^2_{L^2(\Omega)}
		\mathrm{d} t
	\end{displaymath}
	subject to
	\begin{displaymath}
		\left\{
			\begin{array}{ll}
				\frac{\partial v}{\partial t} -\nabla^2 v = u + f,
					& \mathrm{in} \; \Omega \times (0, t_f)\\
				v(\mathbf{x},0)=0, & \mathrm{in} \; \Omega\\
				v(\mathbf{x},t)=0, & \mathrm{on} \; \partial \Omega
					\times(0, t_f),
			\end{array}
		\right.
	\end{displaymath}
	where, for example, $\Omega = (0, 1)^2$, $\beta = 10^{-4}$, and $t_f = 2$.
	The code for the definition  and the solution of the previous problem
	is reported in Figure \ref{heat_control_code}.

	\begin{figure}
	\begin{lstlisting}[language=Python]
	from firedrake import *
	from preconditioner import *
	from control import *

	mesh = UnitSquareMesh(10, 10, 2.0, 2.0)
	space_0 = FunctionSpace(mesh, "Lagrange", 1)

	def forw_diff_operator(trial, test, v, t):
	    return inner(grad(trial), grad(test)) * dx

	def desired_state(test, t):
	    space = test.function_space()
	    mesh = space.mesh()
	    X = SpatialCoordinate(mesh)
	    x = X[0] - 1.0
	    y = X[1] - 1.0

	    v_d = Function(space, name="v_d")
	    v_d.interpolate(t * cos(0.5 * pi * x) cos(0.5 * pi * y))

	    return inner(v_d, test) * dx, v_d

	def force_f(test, t):
	    space = test.function_space()
	    mesh = space.mesh()
	    X = SpatialCoordinate(mesh)
	    x = X[0] - 1.0
	    y = X[1] - 1.0

	    f = Function(space, name="f")
	    f.interpolate(cos(0.5 * pi * x) cos(0.5 * pi * y))

	    return inner(f, test) * dx

	def bc_t(space_0, t):
		return DirichletBC(space_0, 0.0, "on_boundary")

	control_instationary = Control.Instationary(
	    space_0, forw_diff_operator, desired_state=desired_state,
	    force_f=force_f, bcs_v=bc_t, beta=1.0e-4, n_t=10,
	    time_interval=(0.0, 2.0))
	\end{lstlisting}
	\caption{Lines of code required to solve a heat control
		problem.}\label{heat_control_code}
	\end{figure}

The problem is defined in a compact way by providing the weak form representing
	the forward differential operator in space, the boundary conditions on
	the state variable, the desired state, and the force function acting on
	the system. Note that, for stationary control problems, the callables
	defining the desired state and the force function accept as an input
	the test function of the finite element space considered, while the
	callable related to the forward differential operator accepts as inputs
	the trial function, the test function, and the current approximation of
	the state $v$. For instationary problems one has to include also the
	time $t$. Finally the problem is defined by instantiating an
	\texttt{Instationary} object.

Note that the heat control problem defined in Figure \ref{heat_control_code}
	is solved in the time interval $(0,2)$ employing the trapezoidal rule
	discretization over $n_t-1=9$ intervals in time, with a homogeneous initial
	condition for the state and regularization parameter $\beta=10^{-4}$.
	The user can also provide a callable for the definition of
	a different initial condition, passing the argument
	\texttt{initial\_condition}. The discretization in time
	can be set to backward Euler by passing the argument \texttt{CN = False}
	to the call. We would like to note that this example applies implicitly
	the bi-diagonal transformations employed in \cite{Leveque_Pearson_2021}
	for the derivation of the preconditioner, thus resulting in a preconditioner
	structure not currently available, for example, with PETSc field-split
	preconditioning \cite{Brown_Knepley_May_McInnes_Smith}.
	
Control problems which include incompressibility constraints are defined
	by passing to the extra argument \texttt{space\_p}, the space to which
	the pressure belongs. The software assumes that inf--sup stable finite
	element pairs are used.

The linear solvers consist of a Krylov solver preconditioned by suitable
	approximations of the matrices considered. The solvers allow one to
	employ a user-defined preconditioner, adding to the call the argument
	\texttt{P} and passing the solver parameters through the extra argument
	\texttt{solver\_parameters}. The non-linear solver is based on a Picard
	iteration, but it can be set to a Gauss--Newton method by passing the
	argument \texttt{Gauss\_Newton = True} to the definition of the object.
	In the following section, we give details of the in-built preconditioners
	employed in our software.

\subsection{The Preconditioner}\label{sec_3_3}
In this section, we describe the preconditioners implemented in our software
	for solving for the linear systems arising upon discretization of the
	KKT conditions for problems of the type
	\eqref{cost_functional_stationary}--\eqref{stationary_PDE_CO} or
	\eqref{cost_functional_instationary}--\eqref{instationary_PDE_CO}.

When adopting either an optimize-then-discretize or a
	discretize-then optimize approach for solving optimal control problems
	with PDEs as constraints, one has to solve a system of the form
	\eqref{general_block_system}. For complicated problems (like for the control
	of the Navier--Stokes equations), the system in \eqref{general_block_system}
	defines one sub-block of a larger block system. Since the solution of
	systems of the form \eqref{general_block_system} is the most complex task
	in our solver, we describe here only the preconditioners for this type of
	systems. We would like to mention that our preconditioning strategy is
	similar to the field-split approach. In our framework, the
	preconditioner separates the blocks of the discretized system, applying a
	suitable approximation for each ``field''. We would like to
	note that, aside from setting up a field-split approach for the
	preconditioner, our code is deriving the full KKT system given only a few
	callback functions. We refer again to
	\cite{Brown_Knepley_May_McInnes_Smith} for the
	description of the field-split framework, and to
	\cite{Farrell_Kirby_MarchenaMenendez, Farrell_Mitchell_Wechsung,
	Kirby_Mitchell, Rhebergen_Wells_Katz_Wathen} for
	different applications to which it has been applied.

Let $\mathcal{A}$ be given as in \eqref{general_block_system}. Suppose that the
	$(1,1)$-block $A$ is invertible. Then, if we consider the matrix
	\begin{equation}\label{optimal_P}
		\mathcal{P} =
		\left[
			\begin{array}{cc}
				A & 0\\
				B_2 & -S
			\end{array}
		\right],
	\end{equation}
	with invertible (negative) Schur complement $S = C + B_2 A^{-1} B_1^\top$,
	we have that the spectrum of the preconditioned matrix is given by
	$\sigma(\mathcal{P}^{-1}\mathcal{A})=\{ 1 \}$ and the minimal
	polynomial has degree 2, see \cite{Ipsen01, Murphy_Golub_Wathen}.
	From this, an appropriate iterative method should converge in at most
	two iterations in exact arithmetic with this choice of preconditioner.

Despite the optimality of the preconditioner $\mathcal{P}$ in \eqref{optimal_P},
	applying $\mathcal{P}$ is not practical in most of the applications, as
	this requires one to form the Schur complement $S$. In addition, the
	preconditioner requires one to solve exactly for the $(1,1)$-block
	$A$ and for the Schur complement $S$, which may be not invertible
	in some applications and is in general dense even if $\mathcal{A}$ is
	sparse. For these reasons, rather than employing the
	preconditioner $\mathcal{P}$ in \eqref{optimal_P}, one builds
	approximations of the main blocks of $\mathcal{P}$. Specifically,
	given an invertible approximation $\widetilde{A}$ of the
	$(1,1)$-block $A$, one considers an invertible approximation
	$\widetilde{S}$ of the ``perturbed'' Schur complement
	$\widehat{S} = C + B_2 \widetilde{A}^{-1} B_1^\top$. Then,
	the preconditioner that one favours is the following:
	\begin{equation}\label{approx_optimal_P}
		\widetilde{\mathcal{P}} =
		\left[
			\begin{array}{cc}
				\widetilde{A} & 0\\
				B_2 & -\widetilde{S}
			\end{array}
		\right].
	\end{equation}
	In our software, since the preconditioner is non-symmetric,
	the solver of choice is either GMRES \cite{Saad_Schultz} or its flexible
	version \cite{Saad}, restarted after every 10 iterations.

Clearly, the most complex task is finding a suitable approximation of the
	Schur complement $S$. This is because one does not wish to form
	$S$, as mentioned above, but also because the approximation should be
	robust with respect to the problem parameters
	(e.g., the mesh-size, the regularization parameter, the time-step,
	or other parameters defined in the problem). In recent years, the study
	and development of robust preconditioners for systems arising from
	optimal control problems with PDEs as constraints has been a very active
	research area in applied mathematics. We refer the interested reader to
	\cite{Axelsson_Farouq_Neytcheva, Heidel_Wathen,
	Qiu_vanGijzen_vanWingerden_Verhaegen_Vuik, Rees_Dollar_Wathen,
	Schoberl_Zulehner, Stoll_Breiten, Zulehner}.

The preconditioner we adopt in the software is based on the matching
	strategy. This approach was derived for the first time in
	\cite{Pearson_Wathen_2012} for a Poisson control problem, and then
	has been proved to be very effective for the control of stationary
	convection-diffusion problems \cite{Pearson_Wathen_2013}, heat control
	when either employing backward Euler \cite{Pearson_Stoll_Wathen} or
	trapezoidal rule \cite{Leveque_Pearson_2021} in time, and for
	approximating blocks of the sytems arising from Stokes and
	Navier--Stokes control \cite{Leveque_Pearson_2022}. We summarise this
	idea here for the case of a regularized PDE-constrained
	optimization problem. For this class of problems, the linear system
	$\mathcal{A}$ to be solved for is symmetric indefinite,
	that is, $\mathcal{A}$ is defined as in \eqref{general_block_system}
	with $B_1^\top = B_2 = B$, $A \succ 0$, and $C \succ 0$.
	Further, in the first part of the discussion we assume no
	incompressibility constraints are imposed on the state $v$.

With the above assumptions, and writing the Schur complement
	$S = C + B A^{-1} B^\top$ of $\mathcal{A}$, one realises that $S$ is the
	sum of two symmetric positive (semi-) definite matrices, namely,
	$B A^{-1} B^\top$ and $C$. From this expression, the matching strategy
	seeks an approximation of the form
	\begin{equation}\label{hat_S_matching_strategy_general_definition}
		\widehat{S} = (B + \widehat{\Lambda} )
			A^{-1} (B + \widehat{\Lambda})^\top \approx S,
	\end{equation}
	with the matrix $\widehat{\Lambda}$ such that $\widehat{S}$
	`captures' both terms of $S$. As the term $B A^{-1} B^\top$
	already appears in the definition
	\eqref{hat_S_matching_strategy_general_definition} of $\widehat{S}$,
	we wish that the matrix $\widehat{\Lambda}$ is such that
	\begin{equation}\label{lambda_matching_strategy_definition}
		\widehat{\Lambda} A^{-1} \widehat{\Lambda}^\top = C.
	\end{equation}
	It is easy to see that a possible choice for
	\eqref{lambda_matching_strategy_definition} to hold is given by
	\begin{equation}\label{widehat_Lambda}
		\widehat{\Lambda} = C^{\frac{1}{2}} A^{\frac{1}{2}},
	\end{equation}
	where $A^{\frac{1}{2}}$ and $C^{\frac{1}{2}}$ are the (unique) principal
	square roots of $A$ and $C$, respectively. With this choice of
	$\widehat{\Lambda}$, it can be proved that $\frac{1}{2}$ is a lower
	bound for the eigenvalues of the preconditioned Schur complement
	$\widehat{S}^{-1}S$, see \cite[Section 2.11]{Leveque_thesis}.
	Unfortunately, one is not able to prove an upper bound on the spectrum of
	$\widehat{S}^{-1}S$ without further assumptions. Specifically, assuming
	that the ``mixed term''
	$\widehat{\Lambda} A^{-1} B^\top + B A^{-1} \widehat{\Lambda}^\top$
	is at least positive semi-definite, one can prove that $1$ is an upper bound
	for the eigenvalues of $\widehat{S}^{-1}S$, see
	\cite[Section 2.11]{Leveque_thesis}. Note that one can also define
	$\widehat{\Lambda}$ in \eqref{widehat_Lambda} by employing the
	Cholesky decomposition of the matrices $A$ and $C$. In fact, letting
	$A=L_A L_A^\top$ and $C=L_C L_C^\top$, the equality
	\eqref{lambda_matching_strategy_definition} still holds if we choose
	$\widehat{\Lambda} = L_C L_A^\top$. Finally, note that, for
	the optimal control problems with PDEs as constraints
	considered in this work, one does not have to evaluate
	the square root of the matrices $A$ and $C$, as one is the multiple of
	the other. In fact, we have $C = \frac{1}{\beta}A$. Thus,
	we set as matrix $\widehat{\Lambda}= \frac{1}{\sqrt{\beta}} A$. We mention
	that this choice of $\widehat{\Lambda}$ is optimal for all the
	problems with linear PDEs as constraints, without incompressibility
	constraints, considered in this work, as it is possible to prove that the
	mixed term is symmetric positive definite, see
	\cite{Leveque_Pearson_2021, Pearson_Stoll_Wathen,
	Pearson_Wathen_2012, Pearson_Wathen_2013}.

Regarding the $(1,1)$-block, we approximate $A$ either with a fixed number
	of Chebyshev semi-iterations \cite{GolubVargaI, GolubVargaII, Wathen_Rees}
	with a Jacobi splitting, or by employing the Jacobi method.

We conclude this section by describing the approximation of the Schur
	complement for problems with additional incompressibility constraints on
	a state variable $v$ corresponding to velocity. In
	this case, the $(2,1)$-block $B_2$ is a block-diagonal 
	matrix with each diagonal block being the negative divergence matrix,
	with the $(1,2)$-block $B_1^\top$ being the transpose of this matrix,
	while the $(2,2)$-block is the zero matrix. The preconditioner employed
	for the incompressible case is based on the block pressure--convection
	diffusion preconditioner developed in \cite{Leveque_Pearson_2022}.
	Specifically, we employ as a preconditioner
	$\widetilde{\mathcal{P}}$ the block-triangular matrix defined as
	in \eqref{approx_optimal_P}, with the inverse $(1,1)$-block
	$\widetilde{A}$ approximately applied via a fixed number of
	GMRES iterations preconditioned with a block matrix
	derived using the matching strategy above. Finally, the Schur complement
	$S= B_2 A^{-1} B_1^\top$ is approximated with a block-commutator
	argument. Denoting with $\otimes$ the Kronecker product, the Schur
	complement approximation is defined as follows:
	\begin{equation}\label{incompressible_Schur_appr}
		\widetilde{S} = (I_m \otimes K_p) A_p^{-1} (I_m \otimes M_p),
	\end{equation}
	where $K_p$ and $M_p$ are the stiffness and mass matrices on the pressure
	space, and $A_p$ is the discretization on the pressure space of
	the form defined for the $(1,1)$-block $A$ on the velocity space. Here,
	$m\in \mathbb{N}$ is a suitable integer. We would like to mention that
	the approximation $\widetilde{S}$ defined in
	\eqref{incompressible_Schur_appr} can be
	just one of the matrices involved in the whole approximation of the
	Schur complement $S$ (e.g., for time-dependent incompressible problems
	with trapezoidal rule/Crank--Nicolson in time). We refer the reader to
	\cite{Leveque_Pearson_2022} for a more detailed derivation of this
	approximation. Note that, since we have an inner solver for the
	$(1,1)$-block in the case of incompressible problems, we have to employ
	the flexible version of GMRES \cite{Saad} as our outer solver.

\section{Numerical Examples}\label{sec_4}
In this section, we show the applicability of our software to two representative
	(distributed) optimal control problems. The first problem is the control
	of the Poisson equation. This is given as a proof of concept, showing the
	simplicity of the code that our software requires for the definition and
	the solution of the problem. The second is related to the control of
	the time-dependent Navier--Stokes equations with a trapezoidal rule
	discretization in time. This second problem shows the complexity of
	the problem that our software can handle. Regarding this, we report in
	Table \ref{table_problems} the full range of problems for which our software
	has been tested together with the corresponding iterations available.
	Note that our software can solve the optimal control of the Poisson
	equation, the heat equation, the convection--diffusion equation, and the
	Stokes equations, as well as their non-linear variants. We would like to
	note that the our software is not able to deal with a Gauss--Newton
	iteration in the case of solving Navier--Stokes control problems, as the
	block pressure--convection diffusion preconditioner is not well defined
	in this case.

\noindent
\begin{table}[h!]
\begin{center}
	\caption{List of iterative methods that can be employed in our software
		for solving the optimal control problem with the corresponding PDEs.}
	\label{table_problems}
	\renewcommand{\arraystretch}{1.2}
	\begin{footnotesize}
	\begin{tabular}{|c||c||c|}
	\hline
		 & Stationary & Instationary\\
		\hline \hline
		Poisson & Linear, Picard, Gauss--Newton & -- \\
		\hline
		Heat equation & -- & Linear, Picard, Gauss--Newton\\
		\hline
		Convection--diffusion & Linear, Picard, Gauss--Newton &
			Linear, Picard, Gauss--Newton \\
		\hline
		Stokes & Linear, Picard &
			Linear, Picard \\
		\hline
		Navier--Stokes & Picard & Picard \\
		\hline
	\end{tabular}
	\end{footnotesize}
\end{center}
\end{table}

\subsection{Poisson Control}\label{sec_4_1}
We consider here the linear Poisson control problem, with an $L^2$-cost
	functional for the state and the control. Specifically, given the domain
	$\Omega \subset \mathbb{R}^d$ and $\beta > 0$, we consider the following:
	\begin{displaymath}
		\min_{v, u} \| v - v_d \|^2
			+ \dfrac{\beta}{2} \| u \|^2,
	\end{displaymath}
	such that
	\begin{displaymath}
		\left\{
			\begin{array}{ll}
				- \nabla^2 v = u, & \quad \mathrm{in} \; \Omega,\\
				v(x_1, x_2) = 1, & \quad (x_1,x_2) \in \partial \Omega,
			\end{array}
		\right.
	\end{displaymath}
	with desired state given by
	\begin{displaymath}
		v_d(x_1,x_2)= \cos{\left(\frac{\pi x_1}{2}\right)}
			\cos{\left(\frac{\pi x_2}{2}\right)} + 1.
	\end{displaymath}
	For this problem, we set $\Omega = (-1, 1)^2$ and $\beta=10^{-4}$. The
	code for the solution of the above Poisson control problem is given
	in Figure \ref{Poisson_control_code}.

In order to show the efficiency and robustness of the linear solver, we test
	the code on the Poisson control problem defined above. More specifically,
	we run the linear solver for a uniform grid with $2^k+1$ points in
	each spatial dimension, for a range of regularization parameter $\beta$.
	We report the number of GMRES iterations and the CPU times reported
	in seconds in Table \ref{Poisson_control_table}.

As we observe from Table \ref{Poisson_control_table}, the number of iterations
	is roughly constant, with the CPU times scaling linearly with the dimension
	of the problem.

	\begin{figure}
	\begin{lstlisting}[language=Python]
	from firedrake import *
	from preconditioner import *
	from control import *

	mesh = RectangleMesh(10, 10, 2.0, 2.0)
	space_0 = FunctionSpace(mesh, "Lagrange", 1)

	def forw_diff_operator(trial, test, v):
	    return inner(grad(trial), grad(test)) * dx

	def desired_state(test):
	    space = test.function_space()
	    mesh = space.mesh()
	    X = SpatialCoordinate(mesh)
	    x = X[0] - 1.0
	    y = X[1] - 1.0

	    v_d = Function(space, name="v_d")
	    v_d.interpolate(cos(0.5 * pi * x) * cos(0.5 * pi * y) + 1.0)

	    return inner(v_d, test) * dx, v_d

	bc = DirichletBC(space_0, 1.0, "on_boundary")

	control_stationary = Control.Stationary(
	    space_0, forw_diff_operator, desired_state=desired_state,
	    beta=1.0e-4, bcs_v=bc)

	control_stationary.linear_solve()
	\end{lstlisting}
	\caption{Lines of code required to solve a Poisson control
		problem.}\label{Poisson_control_code}
	\end{figure}

	\begin{table}[!htb]
		\caption{Poisson control: number of GMRES iterations and
			elapsed time, for a range of $k$ and
			$\beta$.}\label{Poisson_control_table}
		\begin{footnotesize}
		\begin{center}
			{\begin{tabular}{|c||c|c|c|c|c|c|c|c|c|c|c|c|c|c|}
				\hline
				\multicolumn{1}{|c||}{\phantom{\large{L}}} &
					\multicolumn{2}{c|}{$\beta=10^{0}$} &
					\multicolumn{2}{c|}{$\beta=10^{-1}$} &
					\multicolumn{2}{c|}{$\beta=10^{-2}$} &
					\multicolumn{2}{c|}{$\beta=10^{-3}$} &
					\multicolumn{2}{c|}{$\beta=10^{-4}$} &
					\multicolumn{2}{c|}{$\beta=10^{-5}$} &
					\multicolumn{2}{c|}{$\beta=10^{-6}$}\\
				\cline{2-15}
				$k$ & $\texttt{it}$ & CPU  & $\texttt{it}$ & CPU &
					$\texttt{it}$ & CPU & $\texttt{it}$ & CPU & $\texttt{it}$ &
					CPU & $\texttt{it}$ & CPU & $\texttt{it}$ & CPU \\
				\hline
				\hline
				$~5~$ & 18 & 1.8 & 19 & 1.8 & 18 & 1.9 & 18 & 2.2 & 16 & 2.2 &
					14 & 2.1 & 14 & 2.5 \\
				\hline
				$~6~$ & 19 & 1.4 & 19 & 1.4 & 19 & 1.4 & 18 & 1.4 & 18 & 1.4 &
					16 & 1.2 & 14 & 1.1 \\
				\hline
				$~7~$ & 18 & 3.3 & 19 & 3.5 & 19 & 3.4 & 19 & 3.4 & 19 & 3.6 &
					17 & 3.1 & 16 & 3.0 \\
				\hline
				$~8~$ & 18 & 11 & 20 & 12 & 19 & 11 & 19 & 11 & 19 & 11 &
					17 & 10 & 16 & 10 \\
				\hline
				$~9~$ & 18 & 43 & 20 & 47 & 19 & 45 & 19 & 45 & 20 & 47 &
					 19 & 45 & 17 & 41 \\
				\hline
			\end{tabular}}
		\end{center}
		\end{footnotesize}
	\end{table}

\subsection{Navier--Stokes Control}\label{sec_4_2}
In this section, we consider the time-dependent Navier--Stokes control
	problem defined in \cite[Section\,6.3]{Leveque_Pearson_2022}. This is
	defined as the minimization of the following functional:
	\begin{displaymath}
		\min_{\vec{v}, \vec{u}} \frac{1}{2} \int_0^{t_f} \| \vec{v} - \vec{v}_d
		\|^2 \mathrm{d} t + \frac{\beta}{2} \int_0^{t_f} 
		\| \vec{u} \|^2 \mathrm{d} t
	\end{displaymath}
	subject to
	\begin{displaymath}
		\left\{
			\begin{array}{ll}
				\frac{\partial \vec{v}}{\partial t} - \nu \nabla^2 \vec{v} +
				\vec{v} \cdot \nabla \vec{v} = \vec{u}, & \mathrm{in} \;
					\Omega \times (0, t_f)\\
				\vec{v}(\mathbf{x},0)=\vec{0}, & \mathrm{in} \; \Omega\\
				\vec{v}(\mathbf{x},t)=\vec{g}(\mathbf{x},t), & \mathrm{on} \;
					\partial \Omega \times(0, t_f),
			\end{array}
		\right.
	\end{displaymath}
	where $\Omega = (-1,1)^2$, $\beta=10^{-3}$, $t_f=2$, and the
	boundary conditions are given by
	\begin{displaymath}
		\vec{g}(\mathbf{x},t) =
		\left\{
			\begin{array}{ll}
				\left[t,0\right]^\top & \mathrm{on} \: \partial \Omega_1
					\times (0,1),\\
				\left[1,0\right]^\top & \mathrm{on} \: \partial \Omega_1
					\times [1,2),\\
				\left[0,0\right]^\top & \mathrm{on} \: (\partial \Omega
					\setminus \partial\Omega_1) \times (0,2),
			\end{array}
		\right.
	\end{displaymath}
	with $\partial \Omega_1:= (-1,1) \times \left\{1\right\}$. Here, for
	$\mathbf{x}=(x_1,x_2)$, the desired state is given by
	\begin{displaymath}
		\vec{v}_d(\mathbf{x},t) =
		\left\{
			\begin{array}{ll}
				\vspace{0.5ex}
				c_1 \cos(\frac{\pi t}{2}) \: [(\frac{100}{99})^2 x_2,
					-(\frac{100}{49})^2 (x_1-\frac{1}{2})]^\top
				& \mathrm{if} \: c_1 \geq 0,\\
				\vspace{0.5ex}
				c_2 \cos(\frac{\pi t}{2}) \: [-(\frac{100}{99})^2 x_2,
					(\frac{100}{49})^2 (x_1+\frac{1}{2})]^\top
				& \mathrm{if} \: c_2 \geq 0,\\
				\left[0,0\right]^\top & \mathrm{otherwise},
			\end{array}
		\right.
	\end{displaymath}
	where we set
	\begin{displaymath}
		\begin{array}{c}
			c_1 = 1 - \sqrt{ (\frac{100}{49}( x_1- \frac{1}{2}))^2
				+ (\frac{100}{99}x_2)^2 },\\
			c_2 = 1 - \sqrt{ (\frac{100}{49}( x_1+ \frac{1}{2}))^2
				+ (\frac{100}{99}x_2)^2 }.
		\end{array}
	\end{displaymath}
	We employ a trapezoidal rule for the time discretization, but we could
	just as easily apply backward Euler. The code for the
	solution of this problem is given in Figure \ref{Navier_Stokes_code}.
	
	\begin{figure}
	\begin{lstlisting}[language=Python]
	from firedrake import *
	from preconditioner import *
	from control import *

	import ufl

	mesh = RectangleMesh(10, 10, 2.0, 2.0)

	space_v = VectorFunctionSpace(mesh, "Lagrange", 2)
	space_p = FunctionSpace(mesh, "Lagrange", 1)

	def my_DirichletBC_t_v(space_v, t):
	    if t < 1.0:
	        my_bcs = [DirichletBC(space_v, Constant((t, 0.0)), (4,)),
	                  DirichletBC(space_v, 0.0, (1, 2, 3))]
	    else:
	        my_bcs = [DirichletBC(space_v, Constant((1.0, 0.0)), (4,)),
	                  DirichletBC(space_v, 0.0, (1, 2, 3))]

	    return my_bcs

	def forw_diff_operator_v(trial, test, u, t):
	    nu = 1.0 / 50.0
	    return (nu * inner(grad(trial), grad(test)) * dx
	        + inner(dot(u, grad(trial)), test) * dx)

	def desired_state_v(test, t):
	    space_v = test.function_space()
	    mesh = space_v.mesh()
	    X = SpatialCoordinate(mesh)
	    x = X[0] - 1.0
	    y = X[1] - 1.0

	    a = (100.0 / 49.0) ** 2
	    b = (100.0 / 99.0) ** 2

	    c_1 = 1.0 - sqrt(a * ((x - 0.5) ** 2)
	                     + b * (y ** 2))
	    c_2 = 1.0 - sqrt(a * ((x + 0.5) ** 2)
	                     + b * (y ** 2))
	    v_d = Function(space_v, name="v_d")
	    v_d.interpolate(
	        ufl.conditional(
	            c_1 >= 0.0,
	            c_1 * cos(pi * t / 2.0) * as_vector((b * y,
	                                                 -a * (x - 0.5))),
	            ufl.conditional(
	                c_2 >= 0.0,
	                c_2 * cos(pi * t / 2.0) * as_vector((-b * y,
	                                                     a * (x + 0.5))),
	                as_vector((0.0, 0.0)))),
	    )

	    return inner(v_d, test) * dx, v_d

	control_instationary = Control.Instationary(
	    space_v, forw_diff_operator_v, desired_state=desired_state_v,
	    time_interval=(0.0, 2.0), n_t=10, bcs_v=my_DirichletBC_t_v)

	lambda_v_bounds = (0.3924, 2.0598)
	lambda_p_bounds = (0.5, 2.0)

	control_instationary.incompressible_non_linear_solve(
	    ConstantNullspace(), space_p=space_p,
	    lambda_v_bounds=lambda_v_bounds, lambda_p_bounds=lambda_p_bounds)
	\end{lstlisting}
	\caption{Lines of code required to solve a time-dependent incompressible
		Navier--Stokes control problem, with the trapezoidal rule in
		time.}\label{Navier_Stokes_code}
	\end{figure}
	
The solver employs 20 Chebyshev semi-iterations for applying the approximate
	inverse of mass matrices on the velocity and pressure spaces, and
	2 multigrid cycles for applying the approximate inverse of all
	the other blocks. Note that for incompressible problems, the user has
	to provide as an input the nullspace of the corresponding forward
	stationary Navier--Stokes problem.

We run the solver for different values of the viscosity $\nu$ and regularization
	parameter $\beta$, for the problem defined above. For a given $k$, we
	consider a uniform spatial grid with $2^k+1$ points in each spatial
	direction, and divide the time interval into $2^k - 1$ sub-intervals
	of uniform length. We run the linear solver up to a tolerance of $10^{-6}$
	on the relative residual. We allow for a maximum of 10 Picard
	iterations, and look for a reduction of $10^{-5}$ on the non-linear
	residual. We report the average FGMRES iterations in
	Table \ref{Instationary_NavierStokes_Control}, and the number of non-linear
	iterations in Table \ref{Picard_NavierStokes_Control}.

\begin{table}[!ht]
\caption{Instationary Navier--Stokes control: average FGMRES iterations, for $\nu=\frac{1}{100}$, $\frac{1}{250}$, and $\frac{1}{500}$, and a range of $k$, $\beta$.}\label{Instationary_NavierStokes_Control}
\begin{footnotesize}
\begin{center}
\renewcommand{\arraystretch}{1.2}
{\begin{tabular}{|c||c|c|c||c|c|c||c|c|c|}
\hline
\multicolumn{1}{|c||}{} & \multicolumn{3}{c||}{$\nu=\frac{1}{100}$} & \multicolumn{3}{c||}{$\nu=\frac{1}{250}$} & \multicolumn{3}{c|}{$\nu=\frac{1}{500}$}\\
\cline{2-10}
\multicolumn{1}{|c||}{} & \multicolumn{3}{c||}{$\beta$} & \multicolumn{3}{c||}{$\beta$} & \multicolumn{3}{c|}{$\beta$}\\
\cline{2-10}
$k$ & $10^{-3}$ & $10^{-4}$ & $10^{-5}$ & $10^{-3}$ & $10^{-4}$ & $10^{-5}$ & $10^{-3}$ & $10^{-4}$ & $10^{-5}$ \\
\hline
\hline
$3$ & 11 & 10 & 10 & 10 & 10$\dagger$\footnotemark & 10 & 10 & 10$\dagger$ & 10 \\
\hline
$4$ & 13 & 11 & 10 & 11 & 10 & 9 & 11 & 10 & 9 \\
\hline
$5$ & 18 & 15 & 12 & 16 & 13 & 11 & 14 & 12 & 11 \\
\hline
$6$ & 20 & 17 & 14 & 20 & 16 & 13 & 20 & 14 & 12 \\
\hline
$7$ & 20 & 18 & 17 & 22 & 19 & 16 & 23 & 18 & 14 \\
\hline
\end{tabular}}
\end{center}
\end{footnotesize}
\end{table}

\footnotetext{$\dagger$ means that the Picard iteration did not converge in
	10 iterations. The average number of GMRES iterations and CPU time is
	evaluated over the first 5 Picard iterations.}

\begin{table}[!ht]
\caption{{Instationary Navier--Stokes control: number of Picard iterations required}. In each cell are the Picard iterations for the given $k$, $\nu$, and $\beta=10^{-j}$, $j=3,4,5$.}\label{Picard_NavierStokes_Control}
\begin{small}
\begin{center}
\renewcommand{\arraystretch}{1.2}
\begin{tabular}{|c||ccc|ccc|ccc|}
\hline
$k$ & \multicolumn{3}{c|}{$\nu=\frac{1}{100}$} & \multicolumn{3}{c|}{$\nu=\frac{1}{250}$} & \multicolumn{3}{c|}{$\nu=\frac{1}{500}$} \\
\hline
\hline
$3$ & 8 & 10 & 5 & 9 & $\dagger$ & 5 & 9 & $\dagger$ & 5 \\
\hline
$4$ & 7 & 6 & 6 & 10 & 6 & 7 & 10 & 6 & 7 \\
\hline
$5$ & 5 & 4 & 3 & 5 & 4 & 3 & 6 & 4 & 3 \\
\hline
$6$ & 4 & 4 & 3 & 5 & 4 & 3 & 6 & 4 & 3 \\
\hline
$7$ & 4 & 3 & 3 & 4 & 4 & 4 & 5 & 5 & 4 \\
\hline
\end{tabular}
\end{center}
\end{small}
\end{table}

From Table \ref{Instationary_NavierStokes_Control}, we observe the robustness
	of the preconditioner, obtaining convergence in at most 23 iterations
	for the parameters considered here. Further, from Table
	\ref{Picard_NavierStokes_Control} we observe that the number of
	non-linear iterations decreases as the grid is refined.

\section{Conclusions}\label{sec_5}
In this work we presented an automated software for the numerical solution of
	optimal control problems with a range of PDEs as constraints, when
	no additional algebraic constraints on the variables are imposed.
	The software is based on a Python interface to the Firedrake system.
	This allows the user to provide few lines of a high-level code to define
	the problem considered. The software can solve the optimal control
	of stationary and instationary PDEs, either in the linear or non-linear
	settings. For time-dependent problems, the user can choose to employ
	backward Euler or trapezoidal rule discretization in time. Finally,
	the software can handle problems with additional incompressibility
	constraints on the state variable.

The linear solvers provided in the software are based on advanced
	preconditioning techniques. This results in a drastic speed-up in the
	numerical solution of the problems considered, as in
	\cite{Leveque_Pearson_2021, Leveque_Pearson_2022, Pearson_Stoll_Wathen}
	for instance. The preconditioner is based
	on a field-split strategy, and is applied at the Python level. User-defined
	callables can be passed to the call for different preconditioners, allowing
	for versatility of the software.

Future work includes the development of automated code for the numerical
	solution of optimal control problems with PDEs and additional algebraic
	constraints on the PDE variables. This will consist of devising optimal
	preconditioners for Newton-type methods. Generalizations for the solver
	to handle other classes of PDEs will also be considered.

\section*{Acknowledgements}
The authors gratefully acknowledge Patrick Farrell for a reading of this
	manuscript and for providing very useful suggestions. JWP gratefully
	acknowledges financial support from the Engineering and Physical Sciences
	Research Council (EPSRC) UK grant EP/S027785/1. For the purpose of open
	access, the author has applied a Creative Commons Attribution (CC BY)
	licence to any Author Accepted Manuscript version arising from this
	submission.

\input{references_Control_class}

\end{document}

%% file: references_Control_class.tex
%%%%%%%%%%%%%%%%%%%%%%%% referenc.tex %%%%%%%%%%%%%%%%%%%%%%%%%%%%%%
% sample references
% %
% Use this file as a template for your own input.
%
%%%%%%%%%%%%%%%%%%%%%%%% Springer-Verlag %%%%%%%%%%%%%%%%%%%%%%%%%%
%
% BibTeX users please use
% \bibliographystyle{}
% \bibliography{}
%

%% file: Control_Class.bbl
\begin{thebibliography}{99.}%
% and use \bibitem to create references.
%
% Use the following syntax and markup for your references if 
% the subject of your book is from the field 
% "Mathematics, Physics, Statistics, Computer Science"
%
% Contribution 


% Journal article

\vspace{-1.25em}%!!


\bibitem{Abraham_Behr_Heinkenschloss} F. Abraham, M. Behr, and M. Heinkenschloss: \emph{Shape optimization in steady blood flow: a numerical study of non-Newtonian effects}, Comput. Methods Biomech. Biomed. Engin., 8, pp. 127--137, 2005.

\bibitem{Alnaes_Logg_Olgaard_Rognes_Wells} M. S. Aln\ae{}s, A. Logg, K. B. \O{}lgaard, M. E. Rognes, and G. N. Wells: \emph{Unified form language: A domain-specific language for weak formulations of partial differential equations}, ACM Trans. Math. Softw., 40, Art. 9, 2014.

\bibitem{Arridge} S. R. Arridge: \emph{Optical tomography in medical imaging}, Inverse Problems, 15, pp. R41--R93, 1999.

\bibitem{Axelsson_Farouq_Neytcheva} O. Axelsson, S. Farouq, and M. Neytcheva: \emph{A preconditioner for optimal control problems, constrained by Stokes equation with a time-harmonic control}, J. Comput. Appl. Math., 310, pp. 5--18, 2017.

\bibitem{Barthel_John_Troltzsch} W. Barthel, C. John, and F. Tr\"{o}ltzsch: \emph{Optimal boundary control of a system of reaction diffusion equations}, ZAMM, 90, pp. 966--982, 2010.

\bibitem{Bergounioux_Haddou_Hintermuller_Kunisch} M. Bergounioux, M. Haddou, M. Hinterm\"{u}ller, and K. Kunisch: \emph{A comparison of a Moreau--Yosida-based active set strategy and interior point methods for constrained optimal control problems}, SIAM J. Optim., 11, pp. 495--521, 2000.

\bibitem{Bergounioux_Ito_Kunisch} M. Bergounioux, K. Ito, and K. Kunisch: \emph{Primal-dual strategy for constrained optimal control problems}, SIAM J. Control Optim., 37, pp. 1176--1194, 1999.

\bibitem{Bouchouev_Isakov} I. Bouchouev and V. Isakov: \emph{Uniqueness, stability and numerical methods for the inverse problem that arises in financial markets}, Inverse Problems, 15, pp. R95--R116, 1999.

\bibitem{Brown_Knepley_May_McInnes_Smith} J. Brown, M. G. Knepley, D. A. May, L. C. McInnes, and B. Smith: \emph{Composable linear solvers for multiphysics}, 11th International Symposium on Parallel and Distributed Computing, Munich, Germany, pp. 55--62, 2012.

\bibitem{Casas} E. Casas: \emph{Control of an elliptic problem with pointwise state constraints}, SIAM J. Control Optim., 24, pp. 1309--1318, 1986.

\bibitem{Cheney_Isaacson_Newell} M. Cheney, D. Isaacson, and J. C. Newell: \emph{Electrical impedance tomography}, SIAM Rev., 41, pp. 85--101, 1999.

\bibitem{CHK99} H. Choi, M. Hinze, and K. Kunisch: \emph{Instantaneous control of backward-facing step flows}, Appl. Numer. Math., 31, pp. 133--158, 1999.

\bibitem{Dalcin_Paz_Kler_Cosimo} L. D. Dalcin, R. R. Paz, P. A. Kler, and A. Cosimo: \emph{Parallel distributed computing using Python}, Adv. Water Resour., 34, pp. 1124--1139, 2011.

\bibitem{Egger_Engl} H. Egger and H. W. Engl: \emph{Tikhonov regularization applied to the inverse problem of option pricing: convergence analysis and rates}, Inverse Problems, 21, pp. 1027--1045, 2005.

\bibitem{Ewing_Lin} R. E. Ewing and T. Lin: \emph{A class of parameter estimation techniques for fluid flow in porous media}, Adv. Water Resour., 14, pp. 89--97, 1991.

\bibitem{Farrell_Ham_Funke_Rognes} P. E. Farrell, D. A. Ham, S. W. Funke, and M. E. Rognes: \emph{Automated derivation of the adjoint of high-level transient finite element programs}, SIAM J. Sci. Comput., 35, pp. C369--C393, 2013.

\bibitem{Farrell_Kirby_MarchenaMenendez} P. E. Farrell, R. C. Kirby, and J. Marchena-Men\'{e}ndez: \emph{Irksome: automating Runge--Kutta time-stepping for finite element methods}, ACM Trans. Math. Softw., 47, Art. 30, 2021.

\bibitem{Farrell_Mitchell_Wechsung} P. E. Farrell, L. Mitchell, and F. Wechsung: \emph{An augmented Lagrangian preconditioner for the 3D stationary incompressible Navier--Stokes equations at high Reynolds number}, SIAM J. Sci. Comput., 41, pp. A3073--A3096, 2019.

\bibitem{Funke_Farrell} S. W. Funke and P. E. Farrell: \emph{A framework for automated PDE-constrained optimisation}, arxiv:1302.3894, 2013.

\bibitem{Geiersbach_Wollner} C. Geiersbach and W. Wollner: \emph{A stochastic gradient method with mesh refinement for PDE-constrained optimization under uncertainty}, SIAM J. Sci. Comput., 42, pp. A2750--A2772, 2020.

\bibitem{GolubVargaI} G. H. Golub and R. S. Varga: \emph{Chebyshev semi-iterative methods, successive over-relaxation iterative methods, and second order Richardson iterative methods, Part I}, Numer. Math., 3, pp. 147--156, 1961.

\bibitem{GolubVargaII} G. H. Golub and R. S. Varga: \emph{Chebyshev semi-iterative methods, successive over-relaxation iterative methods, and second order Richardson iterative methods, Part II}, Numer. Math., 3, pp. 157--168, 1961.

\bibitem{Gotschel_Minion} S. G\"{o}tschel and M. L. Minion: \emph{An efficient parallel-in-time method for optimization with parabolic PDEs}, SIAM J. Sci. Comput., 41, pp. C603--C626, 2019.

\bibitem{Griesse_Volkwein} R. Griesse and S. Volkwein: \emph{A primal-dual active set strategy for optimal boundary control of a nonlinear reaction-diffusion system}, SIAM J. Control Optim., 44, pp. 467--494, 2005.

\bibitem{Hasan_Foss_Sagatun} A. Hasan, B. Foss, and S. Sagatun: \emph{Flow control of fluids through porous media}, Appl. Math. Comput., 219, pp. 3323--3335, 2012.

\bibitem{Heidel_Wathen} G. Heidel and A. Wathen: \emph{Preconditioning for boundary control problems in incompressible fluid dynamics}, Numer. Linear Algebra Appl., 26, e2218, 2019.

\bibitem{Hein05} M. Heinkenschloss: \emph{A time-domain decomposition iterative method for the solution of distributed linear quadratic optimal control problems}, J. Comput. Appl. Math., 173, pp. 169--198, 2005.

\bibitem{Hintermuller_Ito_Kunisch} M. Hinterm\"{u}ller, K. Ito, and K. Kunisch: \emph{The primal-dual active set strategy as a semismooth Newton method}, SIAM J. Optim., 13, pp. 865--888, 2002.

\bibitem{Hintermuller_Hinze} M. Hinterm{\"{u}}ller and M. Hinze: \emph{A SQP-semismooth Newton-type algorithm applied to control of the instationary Navier--Stokes system subject to control constraints}, SIAM J. Optim., 16, pp. 1177--1200, 2006.

\bibitem{Hinz00} M. Hinze: \emph{Optimal and instantaneous control of the instationary Navier--Stokes equations}, Habilitation Thesis, Technische Universit\"{a}t Berlin, 2002.

\bibitem{Hinze_Koster_Turek} M. Hinze, M. K\"{o}ster, and S. Turek: \emph{A space--time multigrid method for optimal flow control}. In: G. Leugering, S. Engell, A. Griewank, M. Hinze, R. Rannacher, V. Schulz, M. Ulbrich, S. Ulbrich (eds.), \emph{Constrained optimization and optimal control for partial differential equations}, pp. 147--170, Springer Basel, 2012.

\bibitem{Hinze_Pinnau} M. Hinze and R. Pinnau: \emph{Second-order approach to optimal semiconductor design}, J. Optim. Th. Appl., 133, pp. 179--199, 2007.

\bibitem{Hinze_Pinnau_Ulbrich_Ulbrich} M. Hinze, R. Pinnau, M. Ulbrich, and S. Ulbrich: \emph{Optimization with PDE constraints}, Springer Dordrecht, 2009.

\bibitem{Ipsen01} I. C. F. Ipsen: \emph{A note on preconditioning nonsymmetric matrices}, SIAM J. Sci. Comput., 23, pp. 1050--1051, 2001.

\bibitem{Ito_Kunisch_2003} K. Ito and K. Kunisch: \emph{Semi-smooth Newton methods for state-constrained optimal control problems}, Systems Control Lett., 50, pp. 221--228, 2003.

\bibitem{Kirby_Mitchell} R. C. Kirby and L. Mitchell: \emph{Solver composition across the PDE/linear algebra barrier}, SIAM J. Sci. Comput., 40, pp. C76--C98, 2018.

\bibitem{Klibanov_Lucas} M. V. Klibanov and T. R. Lucas: \emph{Numerical solution of a parabolic inverse problem in optical tomography using experimental data}, SIAM J. Appl. Math., 59, pp. 1763--1789, 1999.

\bibitem{Leveque_thesis} S. Leveque: \emph{Preconditioned iterative methods for optimal control problems with time-dependent PDEs as constraints}, D.Phil. Thesis, University of Edinburgh, 2022.

\bibitem{Leveque_Pearson_2021} S. Leveque and J. W. Pearson: \emph{Fast iterative solver for the optimal control of time-dependent PDEs with Crank--Nicolson discretization in time}, Numer. Linear Algebra Appl., 29, e2419, 2022.

\bibitem{Leveque_Pearson_2022} S. Leveque and J. W. Pearson, \emph{Parameter-robust preconditioning for Oseen iteration applied to stationary and instationary Navier--Stokes control}, SIAM J. Sci. Comput., 44, pp. B694--B722, 2022.

\bibitem{Lions} J. L. Lions: \emph{Optimal control of systems governed by partial differential equations}, Springer Berlin, Heidelberg, 1971.

\bibitem{MaTu02} Y. Maday and G. Turinici: \emph{A parareal in time procedure for the control of partial differential equations}, C. R. Math., 335, pp. 387--392, 2002.

\bibitem{MSS10} T. P. Mathew, M. Sarkis, and C. E. Schaerer: \emph{Analysis of block parareal preconditioners for parabolic optimal control problems}, SIAM J. Sci. Comput., 32, pp. 1180--1200, 2010,

\bibitem{Murphy_Golub_Wathen} M. F. Murphy, G. H. Golub, and A. J. Wathen: \emph{A note on preconditioning for indefinite linear systems}, SIAM J. Sci. Comput., 21, pp. 1969--1972, 2000.

\bibitem{Orozco_Ghattas} C. E. Orozco and O. N. Ghattas: \emph{Massively parallel aerodynamic shape optimization}, Comput. Syst. Eng., 3, pp. 311--320, 1992.

\bibitem{Pearson_Stoll_Wathen} J. W. Pearson, M. Stoll, and A. J. Wathen: \emph{Regularization-robust preconditioners for time-dependent PDE-constrained optimization problems}, SIAM J. Matrix Anal. Appl., 33, pp. 1126--1152, 2012.

\bibitem{Pearson_Wathen_2012} J. W. Pearson and A. J. Wathen: \emph{A new approximation of the Schur complement in preconditioners for PDE-constrained optimization}, Numer. Linear Algebra Appl., 19, pp. 816--829, 2012.

\bibitem{Pearson_Wathen_2013} J. W. Pearson and A. J. Wathen: \emph{Fast iterative solvers for convection--diffusion control problems}, Electron. Trans. Numer. Anal., 40, pp. 294--310, 2013.

\bibitem{Qiu_vanGijzen_vanWingerden_Verhaegen_Vuik} Y. Qiu, M. B. van Gijzen, J.-W. van Wingerden, M. Verhaegen, and C. Vuik: \emph{Preconditioning Navier--Stokes control using multilevel sequentially semiseparable matrix computations}, Numer. Linear Algebra Appl., 28, e2349, 2021.

\bibitem{Quarteroni_Rozza} A. Quarteroni and G. Rozza: \emph{Optimal control and shape optimization of aorto-coronaric bypass anastomoses}, Math. Models Methods Appl. Sci., 13, pp. 1801--1823, 2003.

\bibitem{Rathgeber_Ham_Mitchell_Lange_Luporini_Mcrae_Bercea_Markall_Kelley} F. Rathgeber, D. A. Ham, L. Mitchell, M. Lange, F. Luporini, A. T. T. Mcrae, G.-T. Bercea, G. R. Markall, and P. H. J. Kelly: \emph{Firedrake: automating the finite element method by composing abstractions}, ACM Trans. Math. Softw., 43, Art. 24, 2016.

\bibitem{Rees_Dollar_Wathen} T. Rees, H. S. Dollar, and A. J. Wathen: \emph{Optimal solvers for PDE-constrained optimization}, SIAM J. Sci. Comput., 32, pp. 271--298, 2010.

\bibitem{Rhebergen_Wells_Katz_Wathen} S. Rhebergen, G. N. Wells, R. F. Katz, and A. J. Wathen: \emph{Analysis of block preconditioners for models of coupled magma/mantle dynamics}, SIAM J. Sci. Comput., 36, pp. A1960--A1977, 2014.

\bibitem{Saad} Y. Saad: \emph{A flexible inner--outer preconditioned GMRES algorithm}, SIAM J. Sci. Comput., 14, pp. 461--469, 1993.

\bibitem{Saad_Schultz} Y. Saad and M. H. Schultz: \emph{GMRES: a generalized minimal residual algorithm for solving nonsymmetric linear systems}, SIAM J. Sci. Stat. Comput., 7, pp. 856--869, 1986.

\bibitem{Schiela_Weiser} A. Schiela and M. Weiser: \emph{Superlinear convergence of the control reduced interior point method for PDE constrained optimization}, Comput. Optim. Appl., 39, pp. 369--393, 2008.

\bibitem{Schoberl_Zulehner} J. Sch\"{o}berl and W. Zulehner: \emph{Symmetric indefinite preconditioners for saddle point problems with applications to PDE-constrained optimization problems}, SIAM J. Matrix Anal. Appl., 29, pp. 752--773, 2007.

\bibitem{Stoll_Breiten} M. Stoll and T. Breiten: \emph{A low-rank in time approach to PDE-constrained optimization}, SIAM J. Sci. Comput., 37, pp. B1--B29, 2015.

\bibitem{Sudaryanto_Yortsos} B. Sudaryanto and Y. C. Yortsos: \emph{Optimization of fluid front dynamics in porous media using rate control. I. Equal mobility fluids}, Phys. Fluids, 12, pp. 1656--1670, 2000.

\bibitem{Troltzsch} F. Tr\"{o}ltzsch: \emph{Optimal control of partial differential equations: theory, methods and applications}, Graduate Series in Mathematics, American Mathematical Society, 2010.

\bibitem{Ulbrich_Ulbrich} M. Ulbrich and S. Ulbrich: \emph{Primal-dual interior-point methods for PDE-constrained optimization}, Math. Program., 117, pp. 435--485, 2009.

\bibitem{Wathen_Rees} A. Wathen and T. Rees: \emph{Chebyshev semi-iteration in preconditioning for problems including the mass matrix}, Electron. Trans. Numer. Anal., 34, pp. 125--135, 2009.

\bibitem{Weiser} M. Weiser: \emph{Interior point methods in function space}, SIAM J. Control Optim., 44, pp. 1766--1786, 2005.

\bibitem{Weiser_Schiela} M. Weiser and A. Schiela: \emph{Function space interior point methods for PDE constrained optimization}, PAMM, 4, pp. 43--46, 2004.

\bibitem{Zulehner} W. Zulehner: \emph{Nonstandard norms and robust estimates for saddle point problems}, SIAM J. Matrix Anal. Appl., 32, pp. 536--560, 2011.



\end{thebibliography}
